# Convex optimization in Hilbert space with applications to inverse problems

*Gasnikov A.V. (MIPT, IITP RAS, Moscow),* gasnikov.av@mipt.ru

*Kabanikhin S.I. (ICM&MG, Novosibirsk),* ksi52@mail.ru

*Mohammed A.A.M. (MIPT, Moscow),* a.nafea@science.sohag.edu.eg

*Shishlenin M.A. (ICM&MG, Novosibirsk)* mshishlenin@ngs.ru

**Abstract**

In this paper we propose rather general approaches (based on the gradient descent type methods) to solve convex optimization problems in Hilbert space. We are interested in the case when Hilbert space has infinite dimension. It doesn't typically allow to calculate the exact value of the gradient (Frechet derivative). So that we have a tradeoff between the cost of one iteration and the number of required iterations: one can calculate gradient exactly, so the convergence is fast (we may use fast gradient descent), but the cost of one iteration is large, vice versa, one can calculate gradient roughly so the convergence is slow (we can use only robust simple gradient descent), but the cost of one iteration is cheap. What one should do? In the paper we try to answer for this question. This investigation is motivated by the class of inverse ill-posed problems for elliptic initial-boundary value problems.

**Key words:** convex optimization, inexact oracle, inverse problems.

## 1. Introduction

Following to the works [1, 2, 3] we propose new approaches to solve convex optimization problem in Hilbert space. The main difference from the existing approaches is that we don't approximate infinite-dimensional problem by the finite one (see [1, 3]). We try to solve the problem in Hilbert space (infinite dimensional). But we try to do it with the conception of inexact oracle. That is we use approximation of the problem only when we calculate gradient (Frechet derivative) of the functional or its value. This generates inexactness in gradient calculations. We try to combine all known results in this area and to understand the best way to solve convex optimization problems in Hilbert space with application to inverse problems [2].

The structure of the paper is as follows. In item 2 we described primal approaches (we solve exactly the problem we have) based on contemporary versions of fast gradient descent methods and its adaptive variants. In item 3 we described dual approaches (we solve dual problem) based on the same methods. We try to describe all the methods with the exact estimations of theirs convergence. But every time we have in mind concrete applications. Since that we include in description of algorithms such details that allow method to be more practical.

## 2. Primal approaches

Assume that $q \in H$, where $H$ is a Hilbert space with scalar product denoted by $\langle \,,\, \rangle$ ($H$ isn't necessarily finite). Let's introduce convex functional $J(q)$. In this paper we investigate the following optimization problem

$$J(q) \to \min_q . \tag{1}$$

Let's introduce starting point $y^0$ and

$$R = \left\| y^0 - q_* \right\|_2 ,$$

where $q_*$ is such a solution of (1) that give $R$ the smallest value. We assume that at least one solution exists [1].

Assume that $J(q)$ has Lipchitz Frechet derivative

$$\left\| \nabla J(q_2) - \nabla J(q_1) \right\|_2 \le L \left\| q_2 - q_1 \right\|_2 , \tag{2}$$

where $\|q\|_2^2 = \|q\|_H^2 = \langle q, q \rangle_H$. In (2) we also use that due to the Riesz representation theorem [4], one may considered $\nabla J(q) \in [H \to H]$ to be the element of $H^* = H$.

**Example 1.** Assume that $A: H_1 \to H_2$, $b \in H_2$. Let's consider the following convex optimization problem [1, 2]:

$$J(q) = \frac{1}{2} \left\| Aq - f \right\|_{H_2}^2 \to \min_q .$$

Note that

$$\nabla J(q) = A^* (Aq - f) .$$

Formula (2) is equivalent to

$$\langle Aq, Aq \rangle_{H_2} = \|Aq\|_{H_2}^2 \le L \|q\|_{H_1}^2 = L \langle q, q \rangle_{H_1} , \text{ i.e. } L = \|A\|_{H_1 \to H_2}^2 . \quad \square$$

Now following to [5–7] (most of the ideas below goes back to the pioneer's works of B.T. Polyak, A.S. Nemirovski, Yu.E. Nesterov) we describe optimal (up to absolute constant factor or logarithmic factor in strongly convex case) numerical methods [8] (in terms of the number of ideal calculations of $\nabla J(q)$ and $J(q)$) for solving the problem (1). The rates of convergence that obtained in theorems 1, 2 can be reached (in case of example 1) also by conjugate-gradient methods [1, 2], but we lead these estimates under more general conditions (see, for example, remarks 3, 4).

**Similar Triangular Method** $\mathrm{STM}\left(y^0, L\right)$

**Initialization** ($k=0$)

**Put**

$$A_0 = \alpha_0 = 1/L,\ k=0;\ q^0 = u^0 = y^0 - \alpha_0 \nabla J\left(y^0\right).$$

1. **Put**

$$\alpha_{k+1} = \frac{1}{2L} + \sqrt{\frac{1}{4L^2} + \frac{A_k}{L}},\ A_{k+1} = A_k + \alpha_{k+1},$$

$$y^{k+1} = \frac{\alpha_{k+1} u^k + A_k q^k}{A_{k+1}},\ u^{k+1} = u^k - \alpha_{k+1} \nabla J\left(y^{k+1}\right),\ q^{k+1} = \frac{\alpha_{k+1} u^{k+1} + A_k q^k}{A_{k+1}}.$$

2. If stopping rule doesn't satisfy, put $k := k+1$ and **go to** 1.

If $J\left(q_*\right) = 0$ (see example 1) then stopping rule can have the form $J\left(q^k\right) \le \varepsilon$.

**Theorem 1.** *Assume that (2) holds true. Then*

$$J\left(q^N\right) - J\left(q_*\right) \le \frac{4LR^2}{N^2}.$$

Sometimes it's hardly possible to estimate $L$ that are used in STM. Moreover even when we can estimate $L$ we have to used the worth one (the largest one). Is it possible to change the worth case $L$ to the average one (among all the iterations)? The answer is YES [5] (see ASTM below).

**Adaptive Similar Triangular Method** $\mathrm{ASTM}\left(y^0\right)$

**Initialization** ($k=0$)

**Put**

$$A_0 = \alpha_0 = 1/L_0^0 = 1,\ k=0,\ j_0 = 0;\ q^0 := u^0 := y^0 - \alpha_0 \nabla J\left(y^0\right).$$

**While**

$$\{\ J\left(q^0\right) > J\left(y^0\right) + \left\langle \nabla J\left(y^0\right), q^0 - y^0 \right\rangle + \frac{L_0^{j_0}}{2}\left\|q^0 - y^0\right\|_2^2\ \}$$

**Do**

$$\{\ j_0 := j_0 + 1;\ L_0^{j_0} := 2^{j_0} L_0^0;$$

$$\left(A_0 :=\right)\alpha_0 := \frac{1}{L_0^{j_0}},\ q^0 := u^0 := y^0 - \alpha_0 \nabla J\left(y^0\right)\ \}.$$

1. **Put** $L_{k+1}^0 = L_k^{j_k}/2$, $j_{k+1} = 0$.

$$\alpha_{k+1} := \frac{1}{2L_{k+1}^0} + \sqrt{\frac{1}{4\left(L_{k+1}^0\right)^2} + \frac{A_k}{L_{k+1}^0}},\ A_{k+1} := A_k + \alpha_{k+1},$$

$$y^{k+1} = \frac{\alpha_{k+1} u^k + A_k q^k}{A_{k+1}},\ u^{k+1} = u^k - \alpha_{k+1} \nabla J\left(y^{k+1}\right),\ q^{k+1} = \frac{\alpha_{k+1} u^{k+1} + A_k q^k}{A_{k+1}}.$$

2. **While**

$$\{\ J\left(y^{k+1}\right) + \left\langle \nabla J\left(y^{k+1}\right), q^{k+1} - y^{k+1} \right\rangle + \frac{L_{k+1}^{j_{k+1}}}{2} \left\| q^{k+1} - y^{k+1} \right\|_2^2 < J\left(q^{k+1}\right)\ \}$$

**Do**

$$\{\ j_{k+1} := j_{k+1} + 1;\ L_{k+1}^{j_{k+1}} = 2^{j_{k+1}} L_{k+1}^0;$$

$$\alpha_{k+1} := \frac{1}{2L_{k+1}^{j_{k+1}}} + \sqrt{\frac{1}{4\left(L_{k+1}^{j_{k+1}}\right)^2} + \frac{A_k}{L_{k+1}^{j_{k+1}}}},\ A_{k+1} := A_k + \alpha_{k+1};$$

$$y^{k+1} := \frac{\alpha_{k+1} u^k + A_k q^k}{A_{k+1}},\ u^{k+1} := u^k - \alpha_{k+1} \nabla J\left(y^{k+1}\right),\ q^{k+1} := \frac{\alpha_{k+1} u^{k+1} + A_k q^k}{A_{k+1}}\ \}.$$

3. If stopping rule doesn't satisfy, put $k := k+1$ and **go to** 1.

**Theorem 2.** *Assume that (2) holds true. Then*

$$J\left(q^N\right) - J\left(q_*\right) \le \frac{8LR^2}{N^2}.$$

*The average number of calculations of $J\left(q\right)$ per iteration roughly equals 4 and the average number of calculations Frechet derivative $\nabla J\left(q\right)$ per iteration roughly equals 2.*

**Remark 1 (restarts technique).** Assume that (2) holds true, $J\left(q^*\right) = 0$ and [1]

$$\left\langle \tilde{q}, J''\left(q\right)\tilde{q} \right\rangle \ge \mu \left\langle \tilde{q}, \tilde{q} \right\rangle \text{ for all } \tilde{q}, q \in H\ \ (\mu > 0). \tag{3}$$

For example 1 (3) can be simplified: for all $q \in H$

$$\langle Aq, Aq \rangle_{H_2} = \|Aq\|_{H_2}^2 \geq \mu \|q\|_{H_1}^2 = \mu \langle q, q \rangle_{H_1} \quad (\mu > 0).$$

In this case we may restart the method we choose each time when the value of the goal functional becomes less than one half from starting value at this restart. Since restart criteria is verifiable in practice this construction can be easily implemented in practice. For example, if we choose (A)STM, then by restart construction one can obtain the method that required

$$N = \mathrm{O}\left( \sqrt{\frac{L}{\mu}} \log_2 \left( \frac{\mu R^2}{\varepsilon} \right) \right)$$

calculations of $\nabla J(q)$ (and $J(q)$ in case ASTM) to generate such $q^N$ that guarantee

$$J(q^N) - J(q_*) = J(q^N) \leq \varepsilon .$$

Another approach assumes that we have verifiable stopping criteria $J(\bar{q}_k^{N_k}) - J(q_*) \leq \varepsilon$ and doesn't assume that $J(q_*) = 0$. For example, since

$$J(\bar{q}_k^{N_k}) - J(q_*) \leq \frac{1}{2\mu} \left\| \nabla J(\bar{q}_k^{N_k}) \right\|_2^2 \quad ,$$

we may have such a criteria if we know lower bound on $\mu$. Suppose now that we can estimate $\mu$ from above $\mu \leq \mu_0 \ll L$. Let's consider $\mathrm{RSTM}(y^0, L, \mu_0)$. We introduce the following restart version of this algorithm (here $\mathrm{STM}_q^{N_k}(y_k^0, L) = q^{N_k}$ and $\mathrm{STM}_u^{N_k}(y_k^0, L) = u^{N_k}$):

$$\mathrm{RSTM}(y^0, L, \mu_0)$$

$$k = 0, \; y_0^0 = y^0;$$

**Repeat**

$$\{ \; \bar{q}_k^{N_k} = \frac{1}{2}\mathrm{STM}_q^{N_k}(y_k^0, L) + \frac{1}{2}\mathrm{STM}_u^{N_k}(y_k^0, L), \text{ where } N_k = 2\sqrt{L/\mu_0};$$

$$k := k+1, \; y_{k+1}^0 := \bar{q}_k^{N_k} \; \}$$

**Until**

$$\{ \; J(\bar{q}_k^{N_k}) - J(q_*) \leq \varepsilon \; \}.$$

One can show [9] that after $k = \mathrm{O}\left(\sqrt{\mu_0/\mu}\log_2\left(L\left(J\left(y^0\right)-J\left(q_*\right)\right)/\left(\mu\varepsilon\right)\right)\right)$ restarts $\mathrm{RSTM}\left(y^0, L, \mu_0\right)$ will stop and the total number of calculations of $\nabla J\left(q\right)$ (and $J\left(q\right)$ in case ASTM) can be estimated as follows

$$\mathrm{O}\left(\sqrt{\frac{L}{\mu}}\log_2\left(\frac{L\left(J\left(y^0\right)-J\left(q_*\right)\right)}{\mu\varepsilon}\right)\right).$$

Since (3) is true one has

$$\frac{\mu}{2}\left\|\bar{q}_k^{N_k}-q_*\right\|_2^2 \le J\left(\bar{q}_k^{N_k}\right)-J\left(q_*\right)\le\varepsilon .$$

Therefore,

$$\left\|\bar{q}_k^{N_k}-q_*\right\|_2 \le \sqrt{\frac{2\varepsilon}{\mu}} .$$

Hence after ( $\varepsilon := \mu\varepsilon^2/2$ )

$$\mathrm{O}\left(\sqrt{\frac{L}{\mu}}\log_2\left(\frac{2L\left(J\left(y^0\right)-J\left(q_*\right)\right)}{\mu^2\varepsilon^2}\right)\right)$$

calculations of $\nabla J\left(q\right)$ one can obtain

$$\left\|\bar{q}_k^{N_k}-q_*\right\|_2 \le \varepsilon .$$

For the adaptive case in $\mathrm{RASTM}\left(y^0, \mu_0\right)$ we have to replace $N_k = 2\sqrt{L/\mu_0}$ (since we don't know $L$ ) by $N_k$ – the smallest natural number that satisfies $A_{N_k} \ge 4/\mu_0$ (this result was obtained by D.I. Kamzolov, 2017). All the estimates hold true up to constant factors.

Note that restart techniques (described above) don't assume that we know $\mu > 0$ . If (2), (3) holds true (now there is no need to assume that $J\left(q^*\right)=0$ or something else) and we know $L$, $\mu > 0$ then we can modify $\mathrm{STM}\left(y^0, L\right)$ by changing

$$\alpha_{k+1} = \frac{1+A_k\mu}{2L}+\sqrt{\frac{1+A_k\mu}{4L^2}+\frac{A_k\cdot\left(1+A_k\mu\right)}{L}}$$

and if we know only $\mu > 0$ then we can modify $\mathrm{ASTM}\left(y^0\right)$ by changing

$$\alpha_{k+1} := \frac{1+A_k\mu}{2L_{k+1}^{j_{k+1}}} + \sqrt{\frac{1+A_k\mu}{4\left(L_{k+1}^{j_{k+1}}\right)^2} + \frac{A_k\cdot(1+A_k\mu)}{L_{k+1}^{j_{k+1}}}}.$$

Note that one should also change a bit the structure of methods itself (see details in [5]). For such modifications of $\mathrm{STM}\left(y^0, L\right)$ and $\mathrm{ASTM}\left(y^0\right)$ we have the following generalization of theorems 1, 2

$$J\left(q^N\right) - J\left(q_*\right) \le \min\left\{\frac{4LR^2}{N^2}, LR^2\exp\left(-\frac{N}{2}\sqrt{\frac{\mu}{2L}}\right)\right\}, \text{ // for } \mathrm{STM}\left(y^0, L\right)$$

$$J\left(q^N\right) - J\left(q_*\right) \le \min\left\{\frac{8LR^2}{N^2}, 2LR^2\exp\left(-\frac{N}{2}\sqrt{\frac{\mu}{2L}}\right)\right\}. \text{ // for } \mathrm{ASTM}\left(y^0\right)$$

The last method (as it was in theorem 2) requires in average 4 calculations of $J(q)$ per iteration and 2 calculations of $\nabla J(q)$ per iteration. □

**Remark 2.** In the literature one can typically meet non-accelerated simple gradient descent method $\mathrm{GD}\left(q^0 = y^0, L\right)$ [1, 2, 10]

$$q^{k+1} = q^k - \frac{1}{L}\nabla J\left(q^k\right) \tag{4}$$

or $\overline{\mathrm{GD}}\left(y^0, L\right)$, $q^0 = 0$

$$\begin{cases} y^{k+1} = y^k - \dfrac{1}{L}\nabla J\left(y^k\right), \\ q^{k+1} = \dfrac{k}{k+1}q^k + \dfrac{1}{k+1}y^{k+1}. \end{cases} \tag{5}$$

In the case (2), (3) method (4) requires $\mathrm{O}\left(\left(L/\mu\right)\log_2\left(2LR^2/\left(\mu\varepsilon^2\right)\right)\right)$ calculations of $\nabla J(q)$ for $\left\|q^N - q_*\right\|_2 \le \varepsilon$. In the case (2) method (5) requires $\mathrm{O}\left(LR^2/\varepsilon\right)$ calculations of $\nabla J(q)$ for $J\left(q^N\right) - J\left(q_*\right) \le \varepsilon$. Note that for the (A)STM these quantities are smaller

$$\mathrm{O}\left(\sqrt{\frac{L}{\mu}}\log_2\left(\frac{2LR^2}{\mu\varepsilon^2}\right)\right) \text{ and } \mathrm{O}\left(\sqrt{\frac{LR^2}{\varepsilon}}\right). \tag{6}$$

In reality in (6) for ASTM (analogously for RASTM) one can insert the average $L$ among all the iterations. This could be much smaller than the worth one.

One can easily propose adaptive version of GD and $\overline{\text{GD}}$. For (5) let's describe $\text{AGD}\left(q^0 = y^0\right)$.

**<u>Initialization</u>** ($k = 0$)

**Put**

$$L_0^0 = 1,\ k = 0,\ j_0 = 0;\ q^1 = q^0 - \frac{1}{L_0^0}\nabla J\left(q^0\right).$$

**While**

$$\{\ J\left(q^1\right) > J\left(q^0\right) + \left\langle \nabla J\left(q^0\right), q^1 - q^0 \right\rangle + \frac{L_0^{j_0}}{2}\left\|q^1 - q^0\right\|_2^2\ \}$$

**Do**

$$\{\ j_0 := j_0 + 1;\ L_0^{j_0} := 2^{j_0} L_0^0;$$

$$q^1 = q^0 - \frac{1}{L_0^{j_0}}\nabla J\left(q^0\right)\ \}.$$

1. **Put** $L_{k+1}^0 = L_k^{j_k}/2$, $j_{k+1} = 0$.

$$q^{k+1} = q^k - \frac{1}{L_{k+1}^0}\nabla J\left(q^k\right).$$

2. **While**

$$\{\ J\left(q^k\right) + \left\langle \nabla J\left(q^k\right), q^{k+1} - q^k \right\rangle + \frac{L_{k+1}^{j_{k+1}}}{2}\left\|q^{k+1} - q^k\right\|_2^2 < J\left(q^{k+1}\right)\ \}$$

**Do**

$$\{\ j_{k+1} := j_{k+1} + 1;\ L_{k+1}^{j_{k+1}} = 2^{j_{k+1}} L_{k+1}^0;\ q^{k+1} := q^k - \frac{1}{L_{k+1}^{j_{k+1}}}\nabla J\left(q^k\right)\ \}.$$

3. If stopping rule doesn't satisfy, put $k := k + 1$ and **go to** 1.

For (6) let's describe $\overline{\text{AGD}}\left(y^0\right)$.

**<u>Initialization</u>** ($k = 0$)

**Put**

$$L_0^0 = 1,\ k = 0,\ j_0 = 0;\ q^0 = 0;\ y^1 = y^0 - \frac{1}{L_0^0}\nabla J\left(y^0\right).$$

**While**

$$\{\ J\left(y^1\right) > J\left(y^0\right) + \left\langle \nabla J\left(y^0\right), y^1 - y^0 \right\rangle + \frac{L_0^{j_0}}{2}\left\|y^1 - y^0\right\|_2^2\ \}$$

**Do**

$$\{\ j_0 := j_0 + 1;\ L_0^{j_0} := 2^{j_0} L_0^0;$$

$$y^1 = y^0 - \frac{1}{L_0^{j_0}} \nabla J\left(y^0\right)\ \}.$$

1. **Put** $L_{k+1}^0 = L_k^{j_k}/2$, $j_{k+1} = 0$.

$$y^{k+1} = y^k - \frac{1}{L_{k+1}^0} \nabla J\left(y^k\right),$$

$$q^{k+1} = \frac{k}{k+1} q^k + \frac{1}{k+1} y^{k+1}.$$

2. **While**

$$\{\ J\left(y^k\right) + \left\langle \nabla J\left(y^k\right), y^{k+1} - y^k \right\rangle + \frac{L_{k+1}^{j_{k+1}}}{2}\left\|y^{k+1} - y^k\right\|_2^2 < J\left(y^{k+1}\right)\ \}$$

**Do**

$$\{\ j_{k+1} := j_{k+1} + 1;\ L_{k+1}^{j_{k+1}} = 2^{j_{k+1}} L_{k+1}^0;$$

$$y^{k+1} := y^k - \frac{1}{L_{k+1}^{j_{k+1}}} \nabla J\left(y^k\right),\ q^{k+1} = \frac{k}{k+1} q^k + \frac{1}{k+1} y^{k+1}\ \}.$$

3. If stopping rule doesn't satisfy, put $k := k+1$ and **go to** 1.

The same ("one can insert the average $L$ among all the iterations in estimations $\mathrm{O}\left(\left(L/\mu\right)\log_2\left(LR^2/\varepsilon\right)\right)$, $\mathrm{O}\left(LR^2/\varepsilon\right)$" and "the average number of calculations of $J\left(q\right)$ and $\nabla J\left(q\right)$ per one iteration roughly equals 2") one can say about the rates of convergence for $\mathrm{AGD}\left(q^0 = y^0\right)$ and $\overline{\mathrm{AGD}}\left(y^0\right)$.

Note, that if we change in A(GD), $\overline{\mathrm{AGD}}$ $y^{k+1} = y^k - \frac{1}{L}\nabla J\left(y^k\right)$ by $y^{k+1} = y^k - \alpha_k \nabla J\left(y^k\right)$, where $\alpha_k = \arg\min_{\alpha \geq 0} J\left(y^k - \alpha \nabla J\left(y^k\right)\right)$, then the rates of convergence (in general) don't change up to a constant factors. □

**Remark 3 (inexact oracle).** From remark 2 one may conclude that (A)STM is better than (A)GD ($\overline{\text{(A)GD}}$). In terms of the number of iterations (calculations of $\nabla J(q)$ ($J(q)$)) this is indeed so. But the right criterion is the total number of arithmetic operations (a.o.). Unfortunately, (A)STM is more sensitive to the error in calculation of $\nabla J(q)$ ($J(q)$) than $\overline{\text{(A)GD}}$. Now we plan to say in more details about this issue, following to [10, 11]. First of all, let's denote that if $J(q)$ is convex ($\mu$-strongly convex, $\mu \ge 0$ – see (3)) and (2) holds true then for all $q_1$, $q_2$

$$0 \le \left(\frac{\mu}{2}\|q_2 - q_1\|_2^2 \le\right) J(q_2) - J(q_1) - \langle \nabla J(q_1), q_2 - q_1 \rangle \le \frac{L}{2}\|q_2 - q_1\|_2^2.$$

Assume that for (A)STM at each point $y^{k+1}$ we can observe only such approximate values of $J^\delta\left(y^{k+1}\right)$, $\nabla J^\delta\left(y^{k+1}\right)$ that (in ASTM instead of while{} condition in item 2 one should use the right part of the inequality below with $J\left(q^{k+1}\right) \to J^\delta\left(q^{k+1}\right)$, $\delta \to 2\delta$)

$$0 \le \left(\frac{\mu}{2}\left\|q^{k+1} - y^{k+1}\right\|_2^2 \le\right) J\left(q^{k+1}\right) - J^\delta\left(y^{k+1}\right) - \left\langle \nabla J^\delta\left(y^{k+1}\right), q^{k+1} - y^{k+1} \right\rangle \le \frac{L}{2}\left\|q^{k+1} - y^{k+1}\right\|_2^2 + \delta$$

then (A)STM converges as (constants in $\mathrm{O}(\ )$ is smaller then 5)

$$J\left(q^N\right) - J(q_*) = \mathrm{O}\left(\min\left\{\frac{4LR^2}{N^2}, LR^2 \exp\left(-\frac{N}{2}\sqrt{\frac{\mu}{2L}}\right)\right\}\right) + \mathrm{O}(N\delta).$$

Assume that for $\overline{\text{(A)GD}}$ at each point $y^k$ we can observe only such approximate values of $J^\delta\left(y^k\right)$, $\nabla J^\delta\left(y^k\right)$ that (in $\overline{\text{(A)GD}}$ instead of while{} condition in item 2 one should use the right part of the inequality below with $J\left(y^{k+1}\right) \to J^\delta\left(y^{k+1}\right)$, $\delta \to 2\delta$)

$$0 \le \left(\frac{\mu}{2}\left\|y^{k+1} - y^k\right\|_2^2 \le\right) J\left(y^{k+1}\right) - J^\delta\left(y^k\right) - \left\langle \nabla J^\delta\left(y^k\right), y^{k+1} - y^k \right\rangle \le \frac{L}{2}\left\|y^{k+1} - y^k\right\|_2^2 + \delta$$

then $\overline{\text{(A)GD}}$ converges as (constants in $\mathrm{O}(\ )$ is smaller then 5)

$$J\left(q^N\right) - J(q_*) = \mathrm{O}\left(\min\left\{\frac{LR^2}{N}, LR^2 \exp\left(-N\frac{\mu}{2L}\right)\right\}\right) + \mathrm{O}(\delta).$$

Typically in applications (to inverse problems [2]) we have to calculate conjugate operator $A^*$ for calculation of $\nabla J(q) = A^*(Aq - b)$. These lead us rather often to the initial-boundary value problem for the linear system of partial differential equations. In the most of the cases we can solve this system only numerically by choosing properly small size of the grid $\tau$. So we have $\delta = \mathrm{O}\left(\tau^{p}\right)$, where $p = 1,2,...$ corresponds to the order of the approximation (Chapter 4, [10]). But the cost of calculation approximate values, say, $J^{\delta}\left(y^{k}\right)$, $\nabla J^{\delta}\left(y^{k}\right)$ also depends on $\tau$ like $\mathrm{O}\left(\tau^{-r}\right)$, where $r = 1,2,...$ corresponds to the dimension of the problem (we restrict ourselves here by simple explicit-type scheme). The main problem here is that we can obtain only very rough estimations of the constants in the last two expressions $\mathrm{O}(\ )$. Since that one can propose the following practice-aimed version of the mentioned above algorithms. For the desired accuracy $\varepsilon$ we chose $\delta \sim \varepsilon$ for $\overline{(\mathrm{A})\mathrm{GD}}$ and $\delta \sim \varepsilon\sqrt{\max\left\{\mu R^{2}, \varepsilon\right\}}$ for (A)STM and after that $\tau \simeq C_{GD}\varepsilon^{1/p}$ for $\overline{(\mathrm{A})\mathrm{GD}}$ and $\tau \simeq C_{STM} \cdot \left(\varepsilon\sqrt{\max\left\{\mu R^{2}, \varepsilon\right\}}\right)^{1/p}$ for (A)STM. The constant factors $C_{GD}$, $C_{STM}$ are unknown. Since that it is proper to use restart on this constants. Start with $C_{GD} = 1$ and after

$$N \simeq 4\min\left\{\frac{LR^{2}}{\varepsilon}, \frac{L}{\mu}\ln\left(\frac{LR^{2}}{\varepsilon}\right)\right\}$$

iterations verify $J^{\delta}\left(q^{N}\right) \leq \varepsilon$ or $J\left(q^{N}\right) \leq \varepsilon$ if we can calculate $J(q)$ exactly ($J\left(q_{*}\right) = 0$). If $J^{\delta}\left(q^{N}\right) \leq \varepsilon$ put $C_{GD} := C_{GD}/3$. Analogously for $C_{STM}$ with

$$N \simeq 4\min\left\{\sqrt{\frac{LR^{2}}{\varepsilon}}, \sqrt{\frac{L}{\mu}}\ln\left(\frac{LR^{2}}{\varepsilon}\right)\right\}.$$

For adaptive methods we can use here $L = \max\limits_{k=0,..,N} L_{k}^{j_k}$.

Note, that total number of a.o. for $\overline{(\mathrm{A})\mathrm{GD}}$ can be estimated as follows

$$\mathrm{O}\left(\frac{1}{\varepsilon^{r/p}}\min\left\{\frac{LR^{2}}{\varepsilon}, \frac{L}{\mu}\ln\left(\frac{LR^{2}}{\varepsilon}\right)\right\}\right)$$

and for (A)STM as follows

$$\mathrm{O}\left(\frac{1}{\left(\varepsilon\sqrt{\max\left\{\mu R^{2}, \varepsilon\right\}}\right)^{r/p}}\min\left\{\sqrt{\frac{LR^{2}}{\varepsilon}}, \sqrt{\frac{L}{\mu}}\ln\left(\frac{LR^{2}}{\varepsilon}\right)\right\}\right).$$

It is hardly possible to say what is better from these estimations. As we've already mentioned – we don't know even the right order of the constant in $\mathrm{O}(\ )$. But, roughly, we may expect that $\overline{(\mathrm{A})\mathrm{GD}}$ works better for $r > p$ and (A)STM works better for $r < p$.

Note that in example 1 inexactness in some applications can also arise due to $f$ [2, 12]. The mentioned above theory can be used also in this case.

Numerical experiments (fulfilled by Anastasia Pereberina) show that the described above approaches sometimes work very bad in practice due to the large values of $L$. But in these cases we can use another approach that typically allows us to win at least one or two order in the rate of convergence (these trick was proposed to us by A. Birjukov and A. Chernov). The idea is very simple: start with fixed $L$ (say, $L=1$) and run (non adaptive) algorithm. Then put $L:=2L$ and run algorithm with this parameter (initial/starting point is the same for all the starts). Repeat these restarts $L:=2L$ until we begin to observe stabilization (algorithm with $L$ and $L:=2L$ converges in terms of functional values to the same limit). One can show that the total number of required calculations, that should be growth due to the fact that we don't know real $L$, growth at most for 8 times [13]. But typically we can win much more (than we lose due to restarts) because the methods starts to converges with the value of $L$ that is much smaller than the real one. □

**Remark 4 (regularization).** If instead of (1) we consider regularized problem

$$J^{\mu}(q) = J(q) + \frac{\mu}{2}\|q\|_2^2 \to \min_q \tag{1'}$$

with arbitrary positive $\mu \le \varepsilon/\left(2R^2\right)$ and can find such $q^N$ that

$$J^{\mu}\left(q^N\right) - \min_q J^{\mu}(q) \le \varepsilon/2,$$

then

$$J\left(q^N\right) - \min_q J(q) \le \varepsilon.$$

But the problem (1') is $\mu$-strongly convex.

Note that estimation for GD and AGD $\mathrm{O}\left(LR^2/\varepsilon\right)$ can be obtained (up to a logarithmic factor) from $\mathrm{O}\left(\left(L/\mu\right)\ln\left(2LR^2/\left(\mu\varepsilon^2\right)\right)\right)$ under $\mu \simeq \varepsilon/\left(2R^2\right)$ (see above), analogously for STM and ASTM $\mathrm{O}\left(\sqrt{LR^2/\varepsilon}\right)$ can be obtained (up to a logarithmic factor) from $\mathrm{O}\left(\sqrt{L/\mu}\ln\left(2LR^2/\left(\mu\varepsilon^2\right)\right)\right)$. □

**Remark 5 (non-Euclidian set-up).** All the results mentioned above can be generalized for the convex optimization problem in reflexive Banach space with Lipchitz continuous functional

$$J(q) \to \min_{q\in Q}.$$

Not that the set $Q$ is assumed to be of simple structure that allows to "project" on it efficiently. The proper generalization can be found in [5, 6] (see also [10] for $(\mathrm{A})\mathrm{GD}$ and $\overline{(\mathrm{A})\mathrm{GD}}$ – these methods can be applied also to non convex problems, see [14, 15]). □

**Example 2. (convex optimal control problems).** Let's consider the following optimal control problem ( $q \equiv u(\cdot)$ )

$$J(u(\cdot)) = \int_0^T f^0(t, x(t), u(t))\,dt + \Phi(x(T)) \to \min_{u(\cdot)\in U,\, u(\cdot)\subseteq L_2^m[0,T]},$$

$$\frac{dx}{dt} = f(t, x(t), u(t)),\ x(0) = x^0, \tag{DE}$$

where $U$ is a convex set in $\mathbb{R}^m$ (in terms of remark 5 $U \equiv Q$), all the functions are smooth enough (Chapter 8, [1]), $f(t,x,u)$ is a linear functional of $(x,u)$ with coefficients depend only on $t$ and functional $f^0(t,x,u)$ is convex on $(x,u)$. In this case $J(u(\cdot))$ is convex functional [16]. Due to § 5, Chapter 8 [1][1]

$$\nabla J(u(\cdot)) = \left.\frac{\partial H(t,x,u,\psi)}{\partial u}\right|_{x=x(t,u),\, u=u(t),\, \psi=\psi(t,u)} = H_u\big(t, x(t,u(t)), u(t), \psi(t,u(t))\big),$$

where $H = f^0 + \langle \psi, f\rangle$, $x(t,u)$ is the solution of (DE) and $\psi(t,u)$ is the solution of the conjugate system

$$\frac{d\psi}{dt} = -\frac{\partial H(t,x,u,\psi)}{\partial x},\ \psi(T) = \nabla\Phi(x(T,u)). \tag{CDE}$$

Unfortunately, one can't calculate precisely gradient since one should solve two system of ordinary differential equations (DE), (CDE). But one can solve these two systems by introducing the same lattice (it is significant [3]) in $t$, with the size of each element $\tau$ : $t^{k+1} - t^k \equiv \tau$, for both of the systems (DE), (CDE)

$$\frac{x(t^{k+1}) - x(t^k)}{\tau} = f\big(t^k, x(t^k), u(t^k)\big),\quad x(t^0) = x(0) = x^0,$$

$$\frac{\psi(t^k) - \psi(t^{k+1})}{\tau} = \frac{\partial H}{\partial x}\big(t^{k+1}, x(t^{k+1}), u(t^{k+1}), \psi(t^{k+1})\big),\quad \psi(T) = \nabla\Phi\big(x(t^{T/h})\big).$$

Here we use the standard Euler's scheme [18] with the quality of approximation $\delta \simeq \tau e^{cT}$ (i.e. $\delta \sim \tau$) and the complexity of one iteration is $\sim \tau^{-1}$ (in terms of remark 3 $r = 1$ and $p = 1$). So using ASTM

---

[1] This means that for all small enough $h(\cdot) \in L_2^m[0,T]$ the following holds true

$$J(u(\cdot) + h(\cdot)) - J(u(\cdot)) = \langle \nabla J(u(\cdot)), h(\cdot)\rangle_{L_2^m[0,T]} + \mathrm{O}\left(\|h(\cdot)\|^2_{L_2^m[0,T]}\right) =$$

$$= \int_0^T H_u\big(t, x(t,u(t)), u(t), \psi(t,u(t))\big) h(t)\,dt + \mathrm{O}\left(\int_0^T h(t)^2 dt\right).$$

with proper choice of $\tau \sim \varepsilon^{3/2}$ one can find $\varepsilon$-solution with the total complexity $\sim \varepsilon^{-2}$. The same result (about the total complexity $\sim \varepsilon^{-2}$) is true for AGD. But proper modification of the last method works also with non convex problems (find local extreme) [17].

Note, that if $f^0(t,x,u)$ is also linear functional of $(x,u)$ with coefficients depend only on $t$ then

$$\frac{\partial H(t,x,u,\psi)}{\partial x} \equiv h_0(t) + h_1(t)\psi .$$

If instead of Euler's scheme one can use Runge–Kutta's schemes of order $p \geq 2$ [17, 18], then the complexity of one iteration is still $\sim \tau^{-1}$ ($r = 1$), therefore (see remark 3) we may expect that ASTM will work better than AGD. □

## 3. Dual approaches

Now we concentrate on example 1. The described below approaches goes back to the Yu.E. Nesterov and A.S. Nemirovski (see historical notes in [19, 20]).

Assume that we have to solve the following convex optimization problem

$$g(q) \to \min_{Aq=f}, \tag{7}$$

where $g(q)$ is 1-strongly convex in $H_1$. We build the dual problem

$$\varphi(\lambda) = \max_q \left\{ \langle \lambda, f - Aq \rangle - g(q) \right\} = \langle \lambda, f - Aq(\lambda) \rangle - g(q(\lambda)) \to \min_\lambda . \tag{8}$$

Note, that $\nabla \varphi(\lambda) = f - Aq(\lambda)$.

Let (A)STM with $y^0 = 0$ for the problem (8) generates points $\left\{y^k\right\}_{k=0}^N$, $\left\{u^k\right\}_{k=0}^N$ and $\left\{\lambda^k\right\}_{k=0}^N$ (in (A)STM we denote the last ones by $\left\{q^k\right\}_{k=0}^N$). Put

$$q^N = \sum_{k=0}^N \frac{\alpha_k}{A_N} q\left(y^k\right).$$

Let $q_*$ be the solution of (7) (this solution is unique due to strong convexity of $g(q)$). Then

$$g\left(q^N\right) - g(q_*) \leq \varphi\left(\lambda^N\right) + g\left(q^N\right).$$

The next theorem [19, 20] allows us to calculate the solution of (7) with prescribed precision.

**Theorem 3.** *Assume that we want to solve the problem (7) by passing to the dual problem (8), according to the formulas mentioned above. Let's use (A)STM to solve (8) with the following stopping rule*

$$\varphi\left(\lambda^N\right)+g\left(q^N\right)\le\varepsilon,\ \left\|Aq^N-f\right\|_{H_2}\le\tilde{\varepsilon}.$$

*Then (A)STM stops by making no more than*

$$6\cdot\max\left\{\sqrt{\frac{L\breve{R}^2}{\varepsilon}},\sqrt{\frac{L\breve{R}}{\tilde{\varepsilon}}}\right\} \tag{9}$$

*iterations, where* $L=\left\|A^*\right\|^2_{H_2\to H_1}=\left\|A\right\|^2_{H_1\to H_2}$, $\breve{R}=\left\|\lambda_*\right\|_{H_2}$, $\lambda_*$ *– solution of the problem (8) (if the solution is not unique than we can choose such a solution* $\lambda_*$ *that minimize* $\breve{R}$ *).*

*For ASTM the average number of calculations of* $\varphi\left(\lambda\right)$ *per iteration roughly equals 4 and the average number of calculations Frechet derivative* $\nabla\varphi\left(\lambda\right)=f-Aq\left(\lambda\right)$ *per iteration roughly equals 2.*

**Example 3 (see [19, 21]).** Let's consider the following optimization problem

$$\frac{1}{2}\left\|q\right\|^2_{H_1}\to\min_{Aq=f}.$$

One can built the dual one

$$\min_{Aq=f}\frac{1}{2}\left\|q\right\|^2_{H_1}=\min_q\max_\lambda\left\{\frac{1}{2}\left\|q\right\|^2_{H_1}+\left\langle f-Aq,\lambda\right\rangle\right\}=$$

$$=\max_\lambda\min_q\left\{\frac{1}{2}\left\|q\right\|^2_{H_1}+\left\langle f-Aq,\lambda\right\rangle\right\}=\max_\lambda\left\{\left\langle f,\lambda\right\rangle-\frac{1}{2}\left\|A^*\lambda\right\|^2_{H_1}\right\}. \tag{10}$$

We assume that $Aq=f$ is compatible, hence for the Fredgolm's theorem it's not possible that there exists such a $\lambda$: $A^*\lambda=0$ and $\left\langle b,\lambda\right\rangle>0$.[2] Hence the dual problem is solvable (but the solution isn't necessarily unique). Let's denote $\lambda_*$ to be the solution of the dual problem

$$\varphi\left(\lambda\right)=\frac{1}{2}\left\|A^*\lambda\right\|^2_{H_1}-\left\langle f,\lambda\right\rangle\to\min_\lambda$$

with minimal $H_2$-norm. Let's introduce (from the optimality condition in (10) for $q$): $q\left(\lambda\right)=A^*\lambda$. Using (A)STM for the dual problem one can find (Theorem 3)

---

[2] Indeed, if there exists such $q$ that $Aq=f$ then for all $\lambda$: $\left\langle Aq,\lambda\right\rangle=\left\langle f,\lambda\right\rangle$. Hence, $\left\langle q,A^*\lambda\right\rangle=\left\langle f,\lambda\right\rangle$. Assume that there exists such a $\lambda$, that $A^*\lambda=0$ and $\left\langle f,\lambda\right\rangle>0$. If it is so we observe a contradiction:

$$0=\left\langle q,A^*\lambda\right\rangle=\left\langle f,\lambda\right\rangle>0.$$

$$\left\| Aq^N - f \right\|_{H_1} = \mathrm{O}\left( \frac{L\breve{R}}{N^2} \right), \tag{11}$$

where $L = \left\| A^* \right\|^2_{H_2 \to H_1} = \left\| A \right\|^2_{H_1 \to H_2}$ (as in example 1), $\breve{R} = \left\| \lambda_* \right\|_{H_2}$.

If one will try to solve the primal problem in example 1

$$\frac{1}{2}\left\| Aq - f \right\|^2_{H_2} \to \min_q$$

by (A)STM, one can obtain the following estimate

$$\left\| Aq^N - f \right\|_{H_2} = \mathrm{O}\left( \frac{\sqrt{L}R}{N} \right), \tag{12}$$

where $L = \left\| A \right\|^2_{H_1 \to H_2}$, $R = \left\| q_* \right\|_{H_1}$. Estimate (12) seems worse than (11). But estimate (12) cannot be improving up to a constant factor [8]. There is no contradiction here, since in general $\breve{R}$ can be big enough, i.e. this parameter is uncontrollable. But in real applications we can hope that this (dual) approach lead us to a faster convergence rate (11). □

**Remark 6.** Indeed, all the mentioned above methods (expect (A)GD) are primal-dual ones [19, 20, 21] (if we use their non strongly convex variants). That is for these methods analogues of theorem 3 holds true with proper modification of (9) for (A)$\overline{\text{GD}}$

$$3 \cdot \max\left\{ \frac{L\breve{R}^2}{\varepsilon}, \frac{L\breve{R}}{\tilde{\varepsilon}} \right\}.$$

This means that we can extend remarks 2, 3 for this approach (with the same sensitivity results as in remark 3). Moreover we can also generalize the problem formulation (7) for more general class of the problems [19, 20] (compare this with remark 5)

$$g(q) \to \min_{Aq = f,\, q \in Q}, \tag{7'}$$

where $q$ belongs to reflexive Banach space (with norm $\| \; \|$) and $g(q)$ is 1-strongly convex in norm $\| \; \|$. The dual problem has the form

$$\varphi(\lambda) = \max_{q \in Q}\left\{ \langle \lambda, f - Aq \rangle - g(q) \right\} = \langle \lambda, f - Aq(\lambda) \rangle - g(q(\lambda)) \to \min_\lambda . \; \square \tag{8'}$$

Let's consider another approach to solve problem (7') [13, 19]. This approaches based on the remark 4. We regularize the dual problem (8') (we use $q^N = q(\lambda^N)$ for the solution of (7'))

$$\varphi^\mu(\lambda) = \varphi(\lambda) + \frac{\mu}{2}\left\| \lambda \right\|^2_{H_2} \to \min_\lambda ,$$

where $\mu \simeq \varepsilon / \left(2\breve{R}^2\right)$. Since we don't know $\breve{R}$ (and therefore $\mu$) we may use restart technique on $\mu$ (see remark 1 and 3 and restarts for $C_{GD}$). Let $q_*$ be the solution of (7'). Since ($q(\lambda)$ is determine by (8))

$$g\left(q(\lambda)\right)+\left\langle \lambda, Aq(\lambda)-f\right\rangle \le g\left(q_*\right)$$

we have

$$g\left(q(\lambda)\right)-g\left(q_*\right) \le \|\lambda\|_{H_2} \left\|Aq(\lambda)-f\right\|_{H_2}.$$

The next theorem [13, 19] allows us to calculate the solution of (7') with prescribed precision.

**Theorem 4.** *Assume that we want to solve the problem (7') by passing to the dual problem (8'), according to the formulas mentioned above. Let's use (A)STM to solve (8') with the following stopping rule*

$$\left\|\lambda^N\right\|_{H_2} \left\|Aq\left(\lambda^N\right)-f\right\|_{H_2} \le \varepsilon, \; \left\|Aq\left(\lambda^N\right)-f\right\|_{H_2} \le \tilde{\varepsilon}.$$

*Then (A)STM stops by making no more than*

$$N \simeq 2\sqrt{\frac{L\cdot\left(\varepsilon+2\breve{R}\tilde{\varepsilon}\right)}{\tilde{\varepsilon}^2}} \ln\left(\frac{8L\max\limits_{q_1,q_2\in Q}\left|g\left(q_2\right)-g\left(q_1\right)\right|\cdot\left(\varepsilon+2\breve{R}\tilde{\varepsilon}\right)}{\varepsilon\cdot\tilde{\varepsilon}^2}\right)$$

*iterations.*

*For ASTM the average number of calculations of $\varphi^\mu(\lambda)$ per iteration roughly equals 4 and the average number of calculations Frechet derivative $\nabla\varphi^\mu(\lambda)=f-Aq(\lambda)+\mu\lambda$ per iteration roughly equals 2.*

**Remark 7.** One can extend remarks 2, 3 for this approach too. □

Now let's describe the main motivating example for this paper.

**Example 4 (inverse problem for elliptic initial-boundary value problem).** Let $u$ be the solution of the following problem (P)

$$u_{xx}+u_{yy}=0,\; x,y\in[0,1],$$

$$u_x(0,y)=0,\; y\in[0,1],$$

$$u(1,y)=q(y),\; y\in[0,1],\; q\in L_2[0,1],\; q(0)=q(1)=0,$$

$$u(x,0)=u(x,1)=0,\; x\in[0,1].$$

And corresponding dual problem (D)

$$\psi_{xx} + \psi_{yy} = 0,\ x, y \in [0,1],$$

$$\psi_x(0, y) = \lambda(y),\ y \in [0,1],\ \lambda \in L_2[0,1],\ \lambda(0) = \lambda(1) = 0,$$

$$\psi(1, y) = 0,\ y \in [0,1],$$

$$\psi(x,0) = \psi(x,1) = 0,\ x \in [0,1].$$

Let's introduce Hilbert space

$$H = \{q \in L_2[0,1] : q(0) = q(1) = 0\},$$

linear operator $A : H \to H$ and conjugate operator $A^* : H \to H$ :

$$Aq(\cdot) = u(0,\cdot) \text{ (i.e. for all } y \in [0,1]\ Aq(y) = u(0, y)\text{)}$$

and

$$A^*\lambda(\cdot) = \psi_x(1,\cdot) \text{ (i.e. for all } y \in [0,1]\ A^*\lambda(y) = \psi_x(1, y)\text{)},$$

where $u$ and $\psi$ are corresponding solutions of the problems (P), (D) [2]. To obtain these formulas one may use the general approach, described, for example, in § 7, Chapter 8 [1] (see also Chapter 4 [3]).

Assume that for the problem (P) we don't know real $q_* = q_*(y)$ but we can measure $f(y) = u(0, y)$. The inverse problem [2] can be formulated as: recognize unknown $q_*$ from observable $f$. So we can reduce this problem to the situation described in examples 1, 3

$$J(q) = \frac{1}{2}\|Aq - f\|_H^2 \to \min_{q \in H},$$

$$\frac{1}{2}\|q\|_H^2 \to \min_{Aq = f} \ //\ \varphi(\lambda) = \frac{1}{2}\|A^*\lambda\|_H^2 - \langle f, \lambda\rangle \to \min_{\lambda}.$$

For example 1 one can obtain that $\nabla J(q) = A^*(Aq - f)$. That is, first of all one should solve (P) and calculate $u(0, y) - f(y)$ ($u$ depends on $q$) then one should solve (D) with $\lambda(y) = u(0, y) - f(y)$. At the end, $\nabla J(q) \in [H \to H] = H$ can be interpreted as a function of $y \in [0,1]$:

$$\nabla J(q)(y) = \psi_x(1, y),$$

where $\psi$ is the solution of (D) with $\lambda(y)$ defined above.

For example 3 one can obtain that $\nabla\varphi(\lambda) = f - A\left(A^{*}\lambda\right)$. That is, first of all one should solve (D) and calculate $\psi_x(1, y)$ ($\psi$ depends on $\lambda$) then one should solve (P) with $q(y) = \psi_x(1, y)$. At the end, $\nabla\varphi(\lambda) \in [H \to H] = H$ can be interpreted as a function of $y \in [0,1]$:

$$\nabla\varphi(\lambda)(y) = f(y) - u(0, y),$$

where $u$ is the solution of (P) with $q(y)$ defined above.

Note, that for this example $L = 1$ [2]. □

The research of A.V. Gasnikov was supported Russian Scientific Found project № 14-50-00150. The research of S.I. Kabanikhin and M.A. Shishlenin was supported by the Ministry of Education and Science of the Russian Federation, by SB RAS interdisciplinary grant 14 “Inverse Problems and Applications: Theory, Algorithms, Software”.